\documentstyle[12pt]{article}

\newcommand{\lap}{\mbox{$\bigtriangleup$}}

\newcommand{\be}{\begin{equation}}
\newcommand{\ee}{\end{equation}}
\newcommand{\kernel}{\frac{1}{ |x-y|^{n-\alpha} } }
\newcommand{\ind}{ (n+\alpha)/(n-\alpha) }

\newcommand{\dfrac}[2]{{\mbox{$\displaystyle{\frac{#1}{#2}}$}}}
\newtheorem{mthm}{Theorem}

\newcommand{\ra}{{\mbox{$\rightarrow$}}}
\newtheorem{thm}{Theorem}[section]

\begin{document}
\title{Qualitative Properties of Solutions for
an Integral Equation}
\author{Wenxiong Chen
\thanks{Partially supported by NSF Grant DMS-0072328 }
\hspace{.2in} Congming Li \hspace{.2in} Biao Ou}
\date{}
\maketitle
\begin{quotation}
    {\bf Abstract } Let $n$ be a positive integer and let $ 0 < \alpha < n.$
    In this paper, we continue our study of the integral equation
\be
       u(x) = \int_{R^{n}} \kernel u(y)^{\ind }dy.
\label{0}
\ee

We mainly consider singular solutions in subcritical, critical,
and super critical cases, and obtain qualitative properties, such
as radial symmetry, monotonicity, and upper bounds for the
solutions.


\end{quotation}
\section{Introduction}
\setcounter{equation}{0} Let $R^{n}$ be the $n-$dimensional
Euclidean space, and let $\alpha$ be a constant satisfying $ 0 <
\alpha < n.$ Consider the integral equation
\begin{equation}
      u(x) = \int_{R^{n}} \kernel u(y)^p dy.
\label{eq}
\end{equation}

When $p = \alpha^* := \dfrac{n + \alpha}{n - \alpha}$, it is the
so-called critical case. It arises as an Euler-Lagrange equation
for a functional under a constraint in the context of the
Hardy-Littlewood-Sobolev inequalities. In his elegant paper [L],
Lieb classified the maximizers of the functional, and thus
obtained the best constant in the Hardy-Littlewood-Sobolev
inequalities. He then posed the classification of all the critical
points of the functional -- the solutions of the integral equation
(\ref{eq}) as an open problem.

In our previous paper, we solved this open problem by using the
method of moving planes. We proved that all the solutions of
(\ref{eq}) are radially symmetric and assume the form
\be
c(\frac{t}{t^2 + |x - x_o|^2})^{(n-\alpha)/2}
\label{1}
\ee
with some constant $c = c(n, \alpha)$, and for some $t > 0$ and
$x_o \in R^n$. We also established the equivalence between the
integral equation and the family of well-known semi-linear partial
differential equations
$$(-\lap)^{\frac{\alpha}{2}} u = u^{\frac{n+\alpha}{n-\alpha}} ,
$$
and therefore classified all the solutions of the PDE.

In this paper, we continue to study the integral equation. We
consider subcritical cases $p < \alpha^*$, super critical cases $p
> \alpha^*$, and singular solutions in all cases.

In section 2, we consider subcritical cases. We first prove the
non-existence theorem.

\begin{mthm}
For $p < \alpha^*$, there does not exist regular positive
solutions of (\ref{eq}).
\label{mthm1}
\end{mthm}

Then we consider solutions with one singularity, and use the
method of moving planes to obtain the radial symmetry and
monotonicity of the solutions. The result also applies to singular
solutions in the critical case.

\begin{mthm}
For $p \leq \alpha^*$, if a solution $u$ of (\ref{eq}) has only
one singularity at a point $x^o$, then it must be radially
symmetric about the same point.
\label{mthm2}
\end{mthm}

In section 3, we study singular solutions in the critical case and
obtain an upper bound for the solutions.

\begin{mthm}
Assume that $u(x)$ is a positive solution of (\ref{eq}) with only
one singularity at $x_o$, then there is a constant C, such that
\be u(x) \leq \frac{C}{|x - x_o|^{\frac{n-\alpha}{2}}} .
\label{4.1}
\ee
\label{mthm3}
\end{mthm}

In section 4, we consider super critical cases and provide
examples of non-radially symmetric solutions.

\section{Subcritical Critical Cases}
\medskip

In this section, we prove Theorem 1 and 2.  We first establish the
non-existence theorem for regular solutions.
\begin{thm}
For $p < \alpha^*$, there does not exist regular positive
solutions of (\ref{eq}).
\label{thm2.1}
\end{thm}
Obviously, the integral equation possesses singular solutions. For
instance, one can verify that, for $p > \dfrac{n}{n - \alpha}$,
$$u(x) = \frac{c}{|x|^{\frac{\alpha}{p-1}}} $$
with some appropriate constant c is a singular solution. However,
we can prove the following

\begin{thm}
For $p \leq \alpha^*$, if a solution $u$ of (\ref{eq}) has only
one singularity at a point $x^o$, then it must be radially
symmetric about the same point.
\label{thm2.2}
\end{thm}

The main ingredient of the proofs are the Kelvin type transform
and the method of moving planes.

Assume that $u$ is a solution of integral equation (\ref{eq}). Let
\be v(x) = \frac{1}{|x|^{n-\alpha}} u(\frac{x}{|x|^2})
\label{K}
\ee
 be the Kelvin type transform of $u(x)$. Then it is a straight
forward calculation to verify that $v(x)$ satisfies the equation
\be
v(x) = \int_{R^n} \frac{1}{|x - y|^{n-\alpha}}
|y|^{-(n-\alpha)(\alpha^* -p)} v^p(y) d\,y .
\label{2.1}
\ee

In the following, we will only present the method of moving planes
in a sketchy way. For more details, please see our previous paper
[CLO].
\bigskip

{\bf The Proof of Theorem \ref{thm2.1}.}
\medskip

Assume that $u(x)$ is a positive regular solution of the integral
equation (\ref{eq}). Let $x_1$ and $x_2$ be any two points in
$R^n$. Since the integral equation is invariant under
translations, we may assume that the midpoint $\dfrac{x_1 +
x_2}{2}$ is at the origin. Let $v(x)$ be the Kelvin type transform
as defined in (\ref{K}). Then $v(x)$ has the desired asymptotic
behavior at infinity $$ v(x) \leq \frac{C}{1 + |x|^{n-\alpha}} .$$
Let
$$x_i^* = \frac{x_i}{|x_i|^2} \;\; i = 1, 2$$
be the inversions of $x_i$. Because of the presence of the
singular term $|y|^{-(n-\alpha)(\alpha^* -p)}$ in the integral
equation (\ref{2.1}), similar to what we did in [CLO], we can use
the method of moving planes to show that $v(x)$ must be radially
symmetric about the origin. In particular, $v(x_1^*) = v(x_2^*)$;
and therefore $u(x_1) = u(x_2)$. Since $x_1$ and $x_2$ are any two
points in $R^n$, we conclude that $u$ must be a constant. This is
impossible. Therefore, (\ref{eq}) does not exist any positive
regular solution.
\bigskip

{\bf The Proof of Theorem \ref{thm2.2}.}
\medskip

Without loss of generality, we may assume that $u(x)$ has only one
singularity at point $e = (1, 0, \cdots, 0)$. We show that $u(x)$
is symmetric and monotone decreasing about any line passing
through $e$. Since we do not know any asymptotic behavior of
$u(x)$ at the infinity, we are not able to apply the method of
moving planes on it. To overcome this difficulty, as usual, we
make a Kelvin type transform (\ref{K}) centered at the origin.
Obviously, $v(x)$ still has a singularity at point $e$ and a
possible singularity at the origin. Let $\epsilon$ be any small
positive number, make another Kelvin type transform centered at
$\epsilon e$, i.e. let
$$w_{\epsilon}(x) = \frac{1}{|x|^{n-\alpha}} v(\frac{x}{|x|^2} + \epsilon e)
.$$

Now $w_{\epsilon}(x)$ has a singularity at $e_{\epsilon} :=
\frac{1}{1 - \epsilon} e$ and a possible singularity at the
inversion point of the origin $0^* := -\frac{1}{\epsilon} e$. Now,
we are able to carry on the method of moving planes as we did in
[CLO] to show that $w_{\epsilon}(x)$ is symmetric and monotone
decreasing about the line $\overline{0^* e}$. Since
$$w_{\epsilon}(x) \ra u(x) \,\, \mbox{ as } \epsilon \ra 0 ,$$
We have shown that $u(x)$ is symmetric and monotone decreasing
about the line $\overline{0e}$. Because we can make the Kelvin
type transform centered at any point around $e$, we prove that
$u(x)$ is symmetric and monotone decreasing about any line passing
through $e$. Therefore $u(x)$ must be radially symmetric and
monotone decreasing about the point $e$. This completes the proof
of the theorem.

\section{Critical Case - Singular Solutions}
\medskip

In this section, we consider singular solutions of the integral
equation (\ref{eq}) in the critical case when $p = \alpha^*$. As
one has seen in the previous section,
$$u(x) = \frac{c}{|x|^{\frac{n-\alpha}{2}}} $$ with a suitable
constant is a singular solution. We will show in fact that any
singular solution can not grow faster than this power of $x$.
\begin{thm}
Assume that $u(x)$ is a positive solution of (\ref{eq}) with only
one singularity at $x_o$, then there is a constant C, such that
\be u(x) \leq \frac{C}{|x - x_o|^{\frac{n-\alpha}{2}}} .
\label{4.1}
\ee
\label{thm4.1}
\end{thm}
\medskip

{\bf Proof.}
\medskip

Without loss of generality, we may assume that the solution $u(x)$
has only one singularity at the origin. Then as we have shown in
section 2, $u(x)$ is radially symmetric and monotone decreasing
about the origin. Let $e$ be any point such that $|e| = 1$. Then
by the integral equation (\ref{eq}), we have, for any $r > 0$,
$$ \begin{array}{ll}
u(re) \geq \int_{B_r(0)} \frac{1}{|re - s\omega|^{n-\alpha}}
[u(s)]^{\frac{n+\alpha}{n-\alpha}} s^{n-1} ds d\omega \\ \geq
[u(r)]^{\frac{n+\alpha}{n-\alpha}} \int_0^r \int_{\partial B_1(0)}
\frac{1}{|re - r\omega|^{n-\alpha}} d\omega s^{n-1} ds \\
= [u(r)]^{\frac{n+\alpha}{n-\alpha}} r^{\alpha} \int_0^1
\int_{\partial B_1(0)} \frac{1}{|e - t\omega|^{n-\alpha}} d\omega
t^{n-1} dt \\ = C r^{\alpha} [u(r)]^{\frac{n+\alpha}{n-\alpha}},
\end{array} $$
with some constant $C$. Here, by the radial symmetry of $u$,
$u(re) = u(r)$ for any $e$. It follows that
$$ u(r) \geq \frac{C}{r^{\frac{n-\alpha}{2}}} . $$
This completes the proof of the theorem.

\section{Super Critical Cases}
\medskip

In the super critical case when $p > \alpha^*$, equation
(\ref{eq}) possesses both symmetric solutions and non-symmetric
solutions.

As we mentioned in the previous section,
$$u(x) = \frac{c}{|x|^{\frac{\alpha}{p-1}}} $$
with some appropriate constant c is a singular symmetric solution.

Now we construct a non-radially symmetric solution. Let $x' =
(x_1, \cdots, x_{n-1}),$ let $u(x)$ be a standard solution in
$R^{n-1}$, i.e.
$$u(x') = c (\frac{1}{1 + |x'|^2})^{\frac{n-1-\alpha}{2}} .$$
Then it satisfies
\be
u(x') = \int_{R^{n-1}} \frac{1}{|x' - y'|^{n-1-\alpha}}
[u(y')]^{\frac{n-1+\alpha}{n-1-\alpha}} dy' .
\label{3.1}
\ee
Let $x = (x', x_n)$, and define
$$\tilde{u}(x) = u(x') .$$
Then one can verify that, for some constant $c$,
\be
\tilde{u}(x) =  c \int_{R^n} \frac{1}{|x - y|^{n-\alpha}}
[\tilde{u}(y)]^{\frac{n-1+\alpha}{n-1-\alpha}} dy .
\label{3.2}
\ee
It follows that a constant multiple of $\tilde{u}$ is an
n-dimensional solution of the integral equation in super critical
case, since $\dfrac{n + \alpha}{n - \alpha} < \dfrac{n - 1 +
\alpha}{n - 1 - \alpha} .$ To see (\ref{3.2}), one simply need to
notice from elementary calculus that
$$\int_{-\infty}^{\infty} \frac{1}{|x - y|^{n-\alpha}} dy_n =
\frac{a}{|x' - y'|^{n-1-\alpha}} , $$ with some constant $a$.

\bigskip

\leftline{\bf Addresses and E-mails }
\medskip

{\em Wenxiong Chen}

Department of Mathematics

Yeshiva University

500 W. 185th St.

New York NY 10033

wchen@ymail.yu.edu
\bigskip

{\em Congming Li}

Department of Applied Mathematics

Campus Box 526

University of Colorado at Boulder

Boulder CO 80309

cli@colorado.edu
\bigskip

{\em Biao Ou}

Department of Mathematics

University of Toledo

Toledo OH 43606

bou@math.utoledo.edu

\end{document}